\documentclass{amsart}
\usepackage[pdftex]{graphicx}	
\usepackage{booktabs}
\usepackage{url}

\newcommand{\beqa}{\begin{eqnarray}}
\newcommand{\eeqa}{\end{eqnarray}}
\newcommand{\beqan}{\begin{eqnarray*}}
\newcommand{\eeqan}{\end{eqnarray*}}

\newcommand{\beq}{\begin{equation}}
\newcommand{\eeq}{\end{equation}}
\newcommand{\beqn}{\beqan}
\newcommand{\eeqn}{\eeqan}
\newtheorem{theo}{Theorem}

\newcommand{\bit}{\begin{itemize}}
\newcommand{\eit}{\end{itemize}}

\begin{document}

%
%

\title[Training Neural Networks]
{Training Neural Networks with an algorithm for piecewise linear functions}

\author{Francisco Barahona}

\address{IBM T. J. Watson Research Center, P. O. Box 218\\
Yorktown Heights, NY 10598, USA\\
barahon@us.ibm.com}

\author{Joao Goncalves}

\address{IBM T. J. Watson Research Center, P. O. Box 218\\
Yorktown Heights, NY 10598, USA\\
jpgoncal@us.ibm.com}

\begin{abstract}
We present experiments on training neural networks with an algorithm
that was originally designed as a subgradient method, namely
the Volume Algorithm. We compare
with two of the most frequently used algorithms like Momentum and Adam.
We also compare with the newly proposed algorithm COCOB.
\end{abstract}


\maketitle

\section{Introduction}	
The Volume Algorithm \cite{vol} is a subgradient method designed
for minimizing a convex piecewise linear function, later it was
extended to convex non-differentiable functions in \cite{Bahiense}.
Since its publication in 2000 it has showed better
convergence than the traditional subgradient method \cite{sbg},
and it has been particularly useful for large-scale instances for which other
algorithms have been impractical, see \cite{gunluk}, \cite{escudero}, for instance.
For these
reasons we decided to adapt it for deep learning where high-order optimization
methods are not suitable,
and to
compare it with gradient descent methods.  
Here
we should mention that for instance, the popular ReLu activation function is 
piece-wise linear and
non-differentiable, so a subgradient method is well suited for this case.
We should also mention that in most neural 
networks the loss function
is not convex, then any method designed for minimizing convex functions
becomes a heuristic, and its performance has to be evaluated in an 
experimental fashion. Here we present experiments with 
several different types of neural networks. 
In the final section we present a summary of the results in 
Table~\ref{sample-table}, which shows that the Volume Algorithm is robust
 and compares very well with the others.
To make our results easy to reproduce we restrict our
experiments to several publicly available neural networks, also we make
our code available as open-source in a GitHub repository \cite{BG}.

Stochastic gradient descent \cite{robbins}, \cite{kiefer}, and its relatives like AdaGrad \cite{duchi}, RMSProp \cite{RMS}, Adam \cite{adam}, Momentum
 \cite{momentum} are among the most frequently used algorithms in machine learning;
 see \cite{ruder} for more references. Also recently COCOB has been proposed
 in \cite{Orabona} with promising results. Here we compare
 the Volume algorithm with two of the most commonly used methods,
 namely Momentum and  Adam. We also include COCOB in this comparison.
 We chose Adam according to the survey \cite{ruder} that in Section 4.10
 reads `Kingma et al.  \cite{adam} show that its bias-correction helps Adam slightly outperform RMSprop towards
the end of optimization as gradients become sparser. Insofar, Adam might be the best overall choice'. We chose Momentum because it is widely used.
We also included COCOB because it is a new algorithm that does not require
an initial learning rate, and it has been reported in \cite{Orabona} 
that COCOB has a training performance that is on-par or better than state-of-the-art algorithms.
 
 The remainder of this paper is as follows. In Section 2 we describe the Volume
 Algorithm. Section 3 contains our experiments. In Section 4 we give some final
 remarks.

\section{The Volume Algorithm} 

We start with an example to motivate the choice of the direction in our
method, and then we describe the algorithm.

\subsection{An example}
\begin{figure}[h]
\begin{center}
\includegraphics[width=12cm, height=6.5cm]{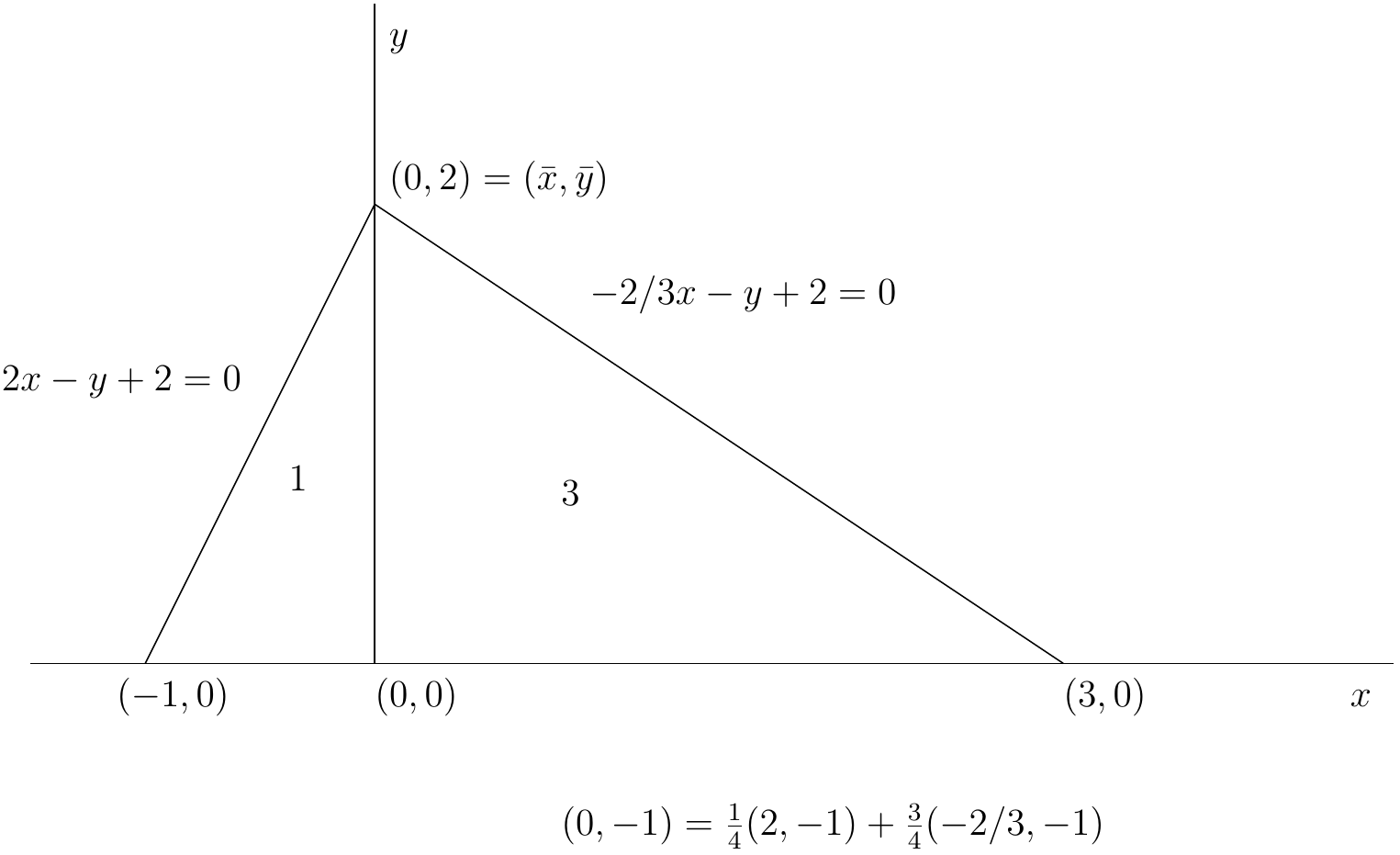}
\end{center}
\caption{Obtaining the direction of improvement as convex combination
of subgradients.}
\label{f1}
\end{figure}

Suppose that we want to maximize the piecewise linear function below
\beqn
f(x,y) = \min \big\{2x -y + 2,
\ - 2/3x -y +2,
\ x + y + 10, 
\ -x + y +10\big\}.
\eeqn
Consider the point $(\bar x, \bar y)=(0,2)$. Here we have
$f(\bar x , \bar y)=0=2 \bar x - \bar y +2= - 2/3 \bar x - \bar y +2$.
Thus at the point $(\bar x, \bar y)$ the function has 
the subgradients $g_1=(2,-1)$ and $g_2=(-2/3, -1)$.
Notice that 
the direction that gives the largest improvement per unit of movement
is $d=(0,-1)$. It can be obtained as a convex combination of
 the subgradients above as 
$d = 1/4 g_1 + 3/4 g_2$. Now we discuss how to obtain the coefficients
of this convex combination. Notice that the area of the triangle with vertices
$(0,2)$, $(0,0)$, $(-1,0)$ is $1$, and 
the area of the triangle with vertices
$(0,2)$, $(0,0)$, $(3,0)$ is $3$. Thus if we use each of these areas
divided by the sum of the areas we obtain the coefficients of the convex
combination, see Figure~\ref{f1}. In higher dimensions we should use volumes
instead of areas. In the next subsection we shall see that this procedure is not
particular to this example, but it corresponds to a method
to obtain the direction of improvement using volumes.

\subsection{Volume and duality}
Now we give a theoretical justification for what we did in the example
above.
Recall that maximizing a concave piecewise linear function can be
written as 
\begin{eqnarray}
&&\hbox{maximize } z \label{lp1}\\
&&\hbox{subject to}\nonumber \\
&&z+ a_i x \le b_i,\hbox{ for } i=1,...,m. \label{lp2}
\end{eqnarray}
The following theorem was published in \cite{vol}.

\begin{figure}[h]
\begin{center}
\includegraphics[width=10cm, height=4.5cm]{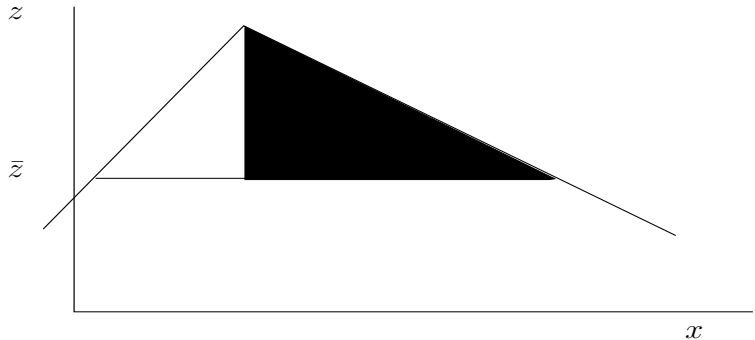}
\end{center}
\caption{Volume and duality.}
\label{f2}
\end{figure}

\begin{theo} \label{th}
Consider the linear program \eqref{lp1}-\eqref{lp2},
let $(\hat z, \hat x)$
be an optimal solution, and suppose that constraints $1, \ldots, m'$,
$m' \le m$, are active at this point. Let $\bar z < \hat z$ and assume that
\beqan
&&z+ a_i x \le b_i,\hbox{ for } i=1,...,m',\\
&& z \ge \bar z
\eeqan
defines a bounded polyhedron. For $1 \le i \le m'$, 
let $\gamma_i$ be the volume between 
the face defined by $z+ a_i x \le b_i$ and the hyperplane defined
by $z = \bar z$. The shaded region in Figure~\ref{f2} illustrates 
such a volume.
Then an optimal dual solution is given by
$$
\lambda_i={\gamma_i \over\sum_{j=1}^{m'} \gamma_j}.
$$
\end{theo}

Roughly speaking, this theorem suggests that given a vector
$(\bar z, \bar x)$ one should look at the
active faces, and compute the volume of their projection over the hyperplane 
$z= \bar z - \epsilon$,
on a neighborhood of $\bar x$, and for some small value of $\epsilon$. 
Let $F_i$ be the face defined by $z+a_i x \le b_i$, and
let $\lambda_i$
be the ratio of the volume
below the face $F_i$ to the total. Then one should compute

$$
(1,0,\cdots,0) - \sum_{i=1}^{m'} \lambda_i (1, a_i).
$$

\begin{figure}[h]
\begin{center}
\includegraphics[width=10cm, height=3.cm]{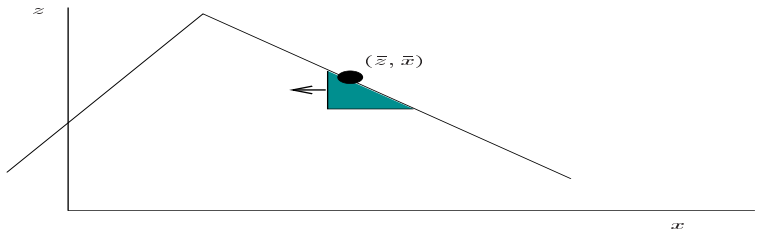}
\end{center}
\caption{Choice of a direction.}
\label{f3}
\end{figure}

If this is the vector of all zeroes, we have a proof of optimality, otherwise we have a direction of improvement,
see Figure~\ref{f3}. This is similar to what we did in the previous
example.
In the next subsection we present an algorithm that
computes approximations of these volumes.

\subsection{The Algorithm}

As previously suggested, for a vector $(\bar z, \bar x)$,
 we could try a set of random points in the hyperplane
 defined by $z=\bar z -\epsilon$. For each of these points
 we can identify the face above it.
Theorem~\ref{th} suggests that we should look at the coefficients 
$\{\lambda_i \}$ as the probabilities that the different faces
will appear. Thus we look at the coefficients 
$\{\lambda_i \}$ as the probabilities 
that a subgradient method will produce
the different subgradients . Then to compute the 
direction of improvement we propose the keep a direction $d$ and 
update it in every iteration as
\beq
d \leftarrow \alpha g + (1-\alpha) d,
\label{dir}
\eeq
where $g$ is the last subgradient produced, and $\alpha$ is a small value.
Therefore if $g_0, \ldots,g_k$ is the sequence of subgradients produced, 
the direction is
\beqn
d =
(1-\alpha)^k g_0 + \alpha (1-\alpha)^{k-1} g_1+ \ldots + 
\alpha(1-\alpha) g_{k-1} + 
\alpha g_k.
\eeqn
The coefficients of this convex combination are an approximation of
the values $\{\lambda_i\}$ of Theorem~\ref{th}. This update scheme
is similar to the one used in Adam \cite{adam}, the difference with our
algorithm is in the choice of the step-size. As we shall see in the next 
section, this makes a significant difference in terms of convergence.

Once the direction $d$ has been chosen, ideally we would do a line search,
but since this is prohibitively expensive, we use the following scheme.
Suppose we use a step size $l$, and let $\bar g$ be the new subgradient
obtained. If the inner product $d \cdot \bar g$ is positive, this means that
the step size was too small, a larger step size would have resulted on
a better solution. We call this a ``green'' iteration.
If $d \cdot \bar g$ is negative, this means that a smaller step size would have
given a better solution. We call this a ``yellow'' iteration. 
Then we keep a ``green-yellow'' indicator $gy$ that is a number
between $0$ and $1$. At each green iteration we update $gy$
as $gy \leftarrow (1-\alpha) gy + \alpha \,1$. At each yellow iteration
we update $gy$ as $gy \leftarrow (1-\alpha) gy + \alpha \,0$.
Each time that $gy > 0.9$ we increase the step size. When $gy < 0.8$
we decrease the step size. We see this as low cost approximation 
of a line search.

Now we can describe the algorithm, it was published first in
\cite{vol}, and a proof of convergence was given in \cite{Bahiense}.

\vskip 0.3cm
\centerline{\bf Volume Algorithm}
\bit
\item Input: 
We assume that we have to maximize a function
$f$ of variables $v$. We also receive 
a value $\alpha \in [0, 1]$, and an initial step-size $s$.
We set upper an upper bound for the step-size $U=2s$, and 
and a lower bound $L=0.2s$. Initially we set $gy=0.5$. 

At each iteration we keep a vector
$d$ that corresponds to a direction. 
We initialize it to be equal to the first subgradient.
At each iteration we receive a new subgradient
$g$ that is used to update the direction. The different steps are below.
\item {\bf Step 1.} Compute $p=d \cdot g$. This is the inner product of the direction
and the new subgradient. 
\newline
If $p > 0$ increase $gy$ as
\hskip 1cm$ gy \leftarrow (1-\alpha) gy + \alpha \, 1,$
\newline
otherwise decrease $gy$ as
\hskip 0.7cm$ gy \leftarrow (1-\alpha) gy + \alpha \, 0.$
\newline
If $gy > 0.9$ increase $s$ as $s \leftarrow 1.01 s$. 
\newline
If $gy < 0.8$ decrease $s$ as $s \leftarrow 0.99 s$.
\newline
If $s > U$, set $s=U$. If $s < L$, set $s=L$.
\newline
Update the direction as
$$ d \leftarrow (1-\alpha) d + \alpha \, g.$$

\item {\bf Step 2.} Update the variables $v$ as
$$ v \leftarrow v + \frac{s}{|| d ||} d.$$
Go to Step 1.
\eit

This is for maximization,
if we want to minimize we update the variables as
$ v \leftarrow v - \frac{s}{|| d ||} d.$ In our experiments
we used $\alpha = 0.1$.

\section{Experiments}

Here we compare the Volume Algorithm with three of the most frequently used
methods, namely Momentum \cite{momentum}, Adam \cite{adam} and
RMSProp \cite{RMS},
we also compare with the new algorithm COCOB \cite{Orabona}.
For the first three algorithms we used their implementation found in 
\cite{TF}, for COCOB we used the implementation in \cite{cocob}.
The Volume Algorithm 
was implemented in TensorFlow. For the choice of the 
learning rate,
for each algorithm
we looked for the best value in 
$\{0.00001, 0.000025, 0.00005, 0.000075, 0.0001, 0.00025,$
    $0.0005, 0.00075, 0.001, 0.0025,$ $0.005,
    0.0075, 0.01,
    0.025,$
    $0.05, 0.075, 0.1\}$. This set of values was used in \cite{Orabona} for their
    comparisons. Recall that COCOB has the interesting feature of not requiring an initial learning rate.
    In the figures below we show the best learning rate used
    for each algorithm.
For Adam we used $\beta_1=0.9$ and $\beta_2=0.999$.    
For Momentum we used $0.9$ for the momentum value.

    We used seven test sets available in the public
    domain, as follows. We have downloaded tensorflow implementations of
    different neural networks in the public domain. We used the same
    number of steps, epochs, batch-size, etc.,  as in these implementations.
    We make our code publicly available in \cite{BG}, and thus these experiments are
    easy to reproduce.
    
    \subsection{MNIST}
\begin{figure}[h!]
\begin{center}
\includegraphics[width=0.5\linewidth, height=0.53\linewidth] {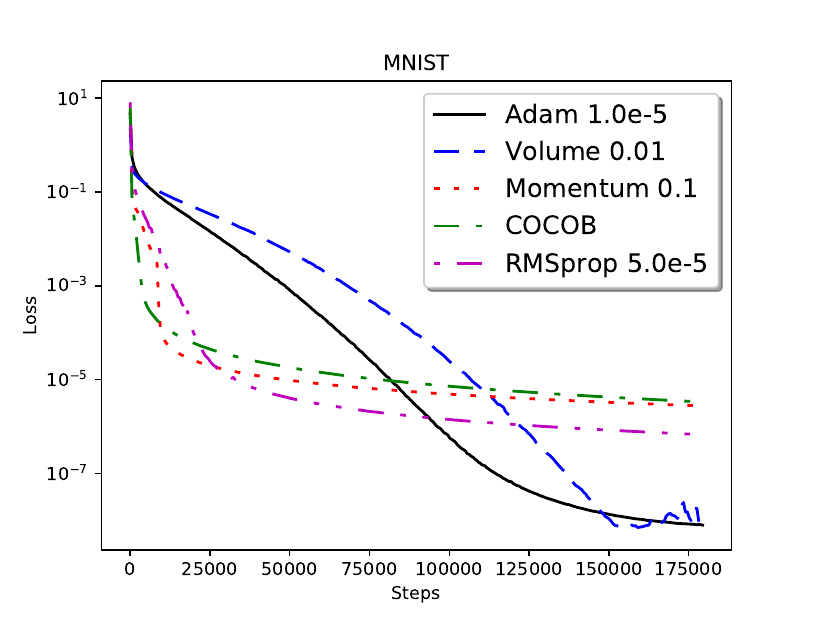}\includegraphics[width=0.5\linewidth, height=0.53\linewidth]{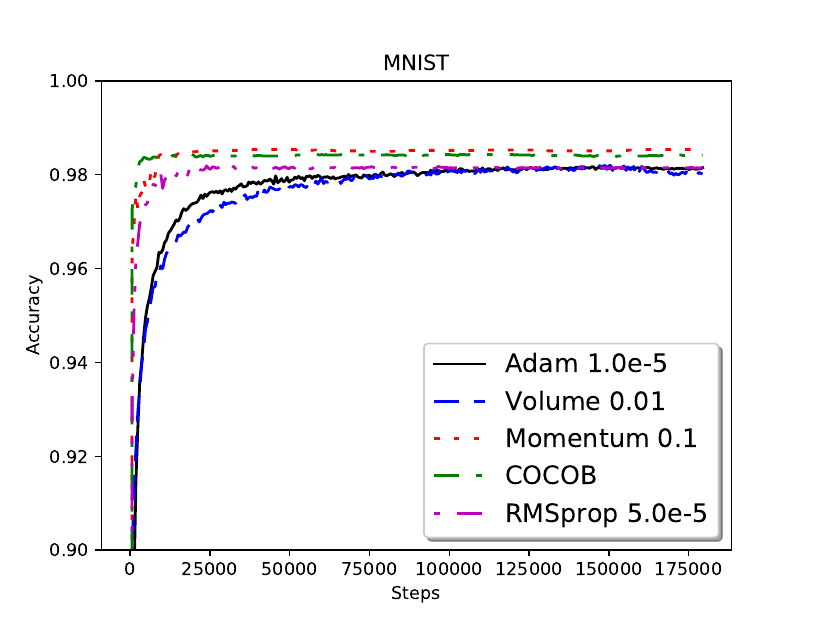}
\end{center}
\caption{MNIST.}
\label{f4}
\end{figure}

    This is a well known digits-recognition dataset, see \cite{mnist}. We used the
    neural network implemented 
    in \cite{cocob}.
    In Figure~\ref{f4} we plot train loss and accuracy
    as functions of the number of steps.

\subsection{CIFAR10}
\begin{figure}[h!]
\begin{center}
\includegraphics[width=0.5\linewidth, height=0.53\linewidth]{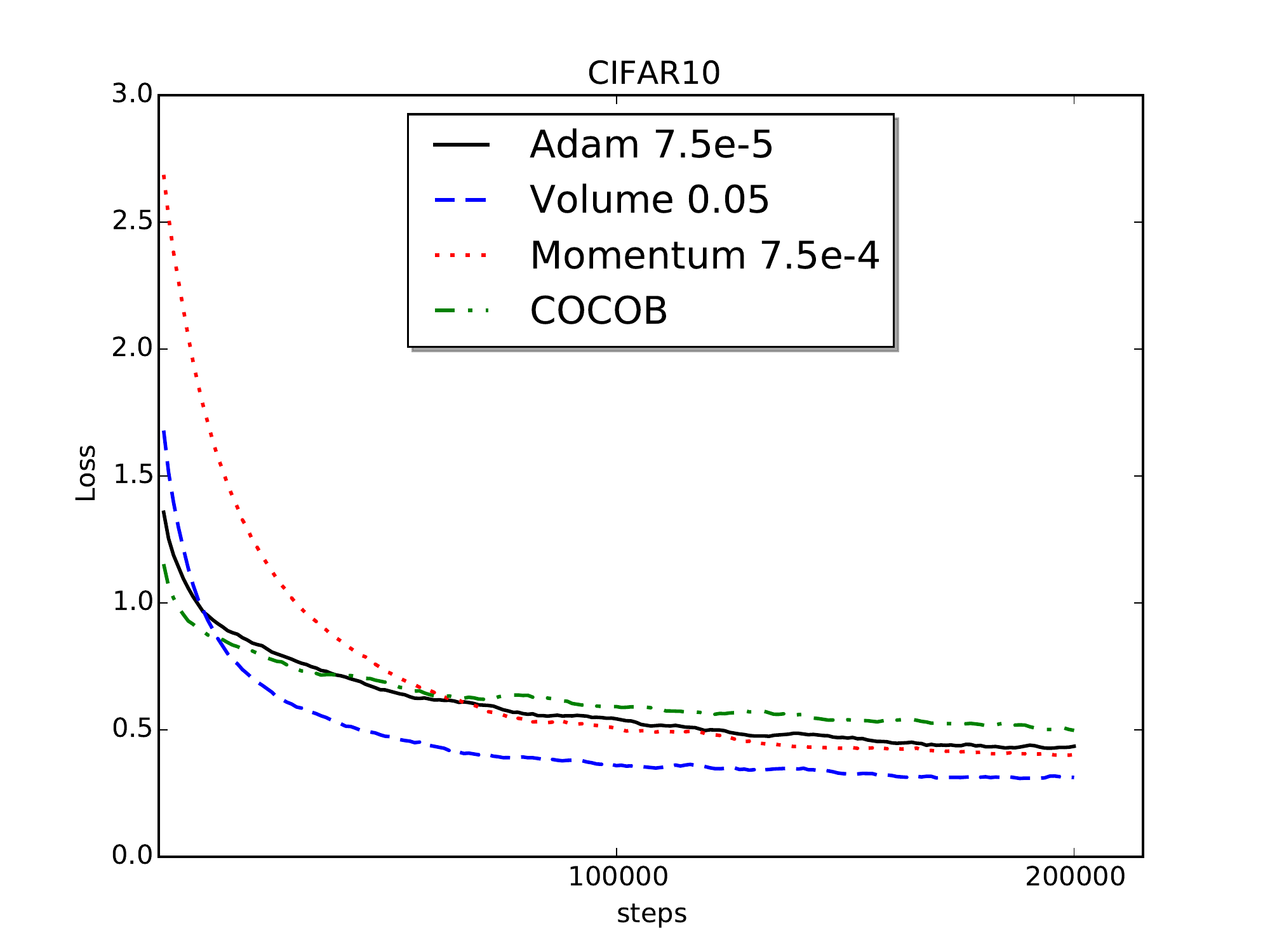}\includegraphics[width=0.5\linewidth, height=0.53\linewidth]{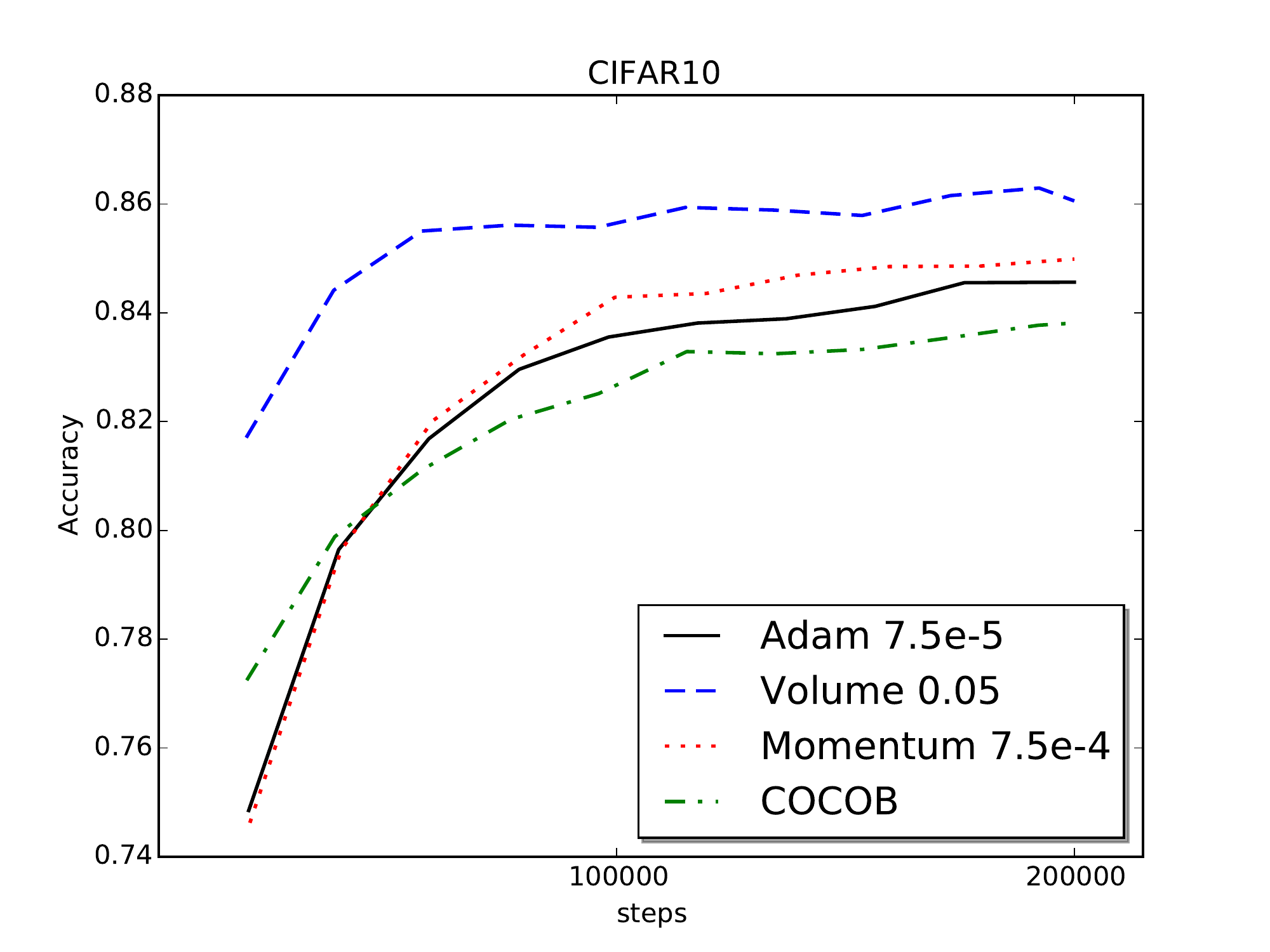}
\end{center}
\caption{CIFAR10.}
\label{f5}
\end{figure}

This is a dataset for classifying 32x32 RGB
images across 10 object categories, see \cite{cifar10}. Here we used the
code from the TensorFlow tutorial \cite{TFcifar}. 
In Figure~\ref{f5} we plot
train loss and accuracy as functions of the number of steps.

\subsection{Recurrent Neural Networks}
\begin{figure}[h]
\begin{center}
\includegraphics[width=0.5\linewidth, height=0.5\linewidth]{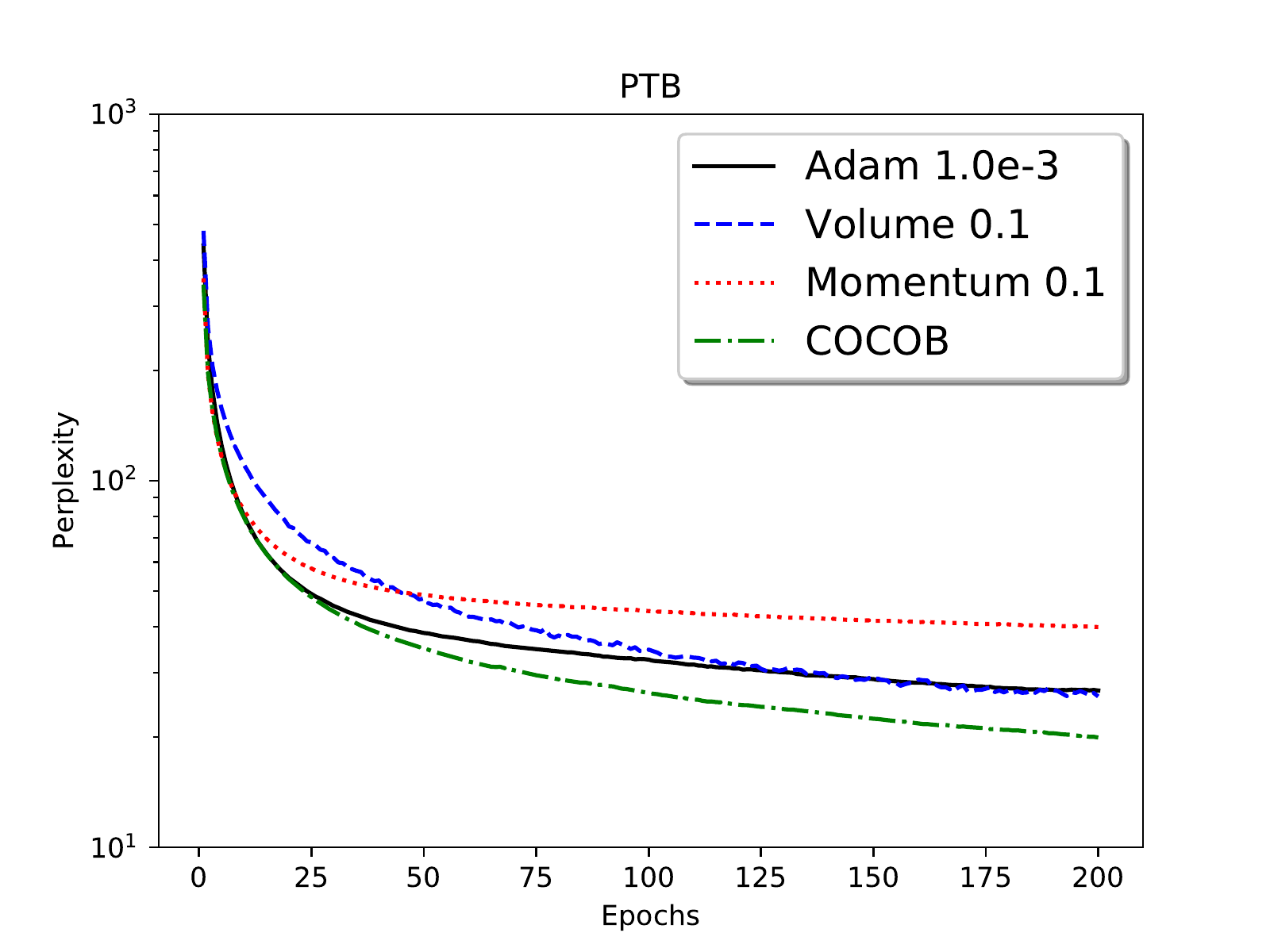}\includegraphics[width=0.5\linewidth, height=0.5\linewidth]{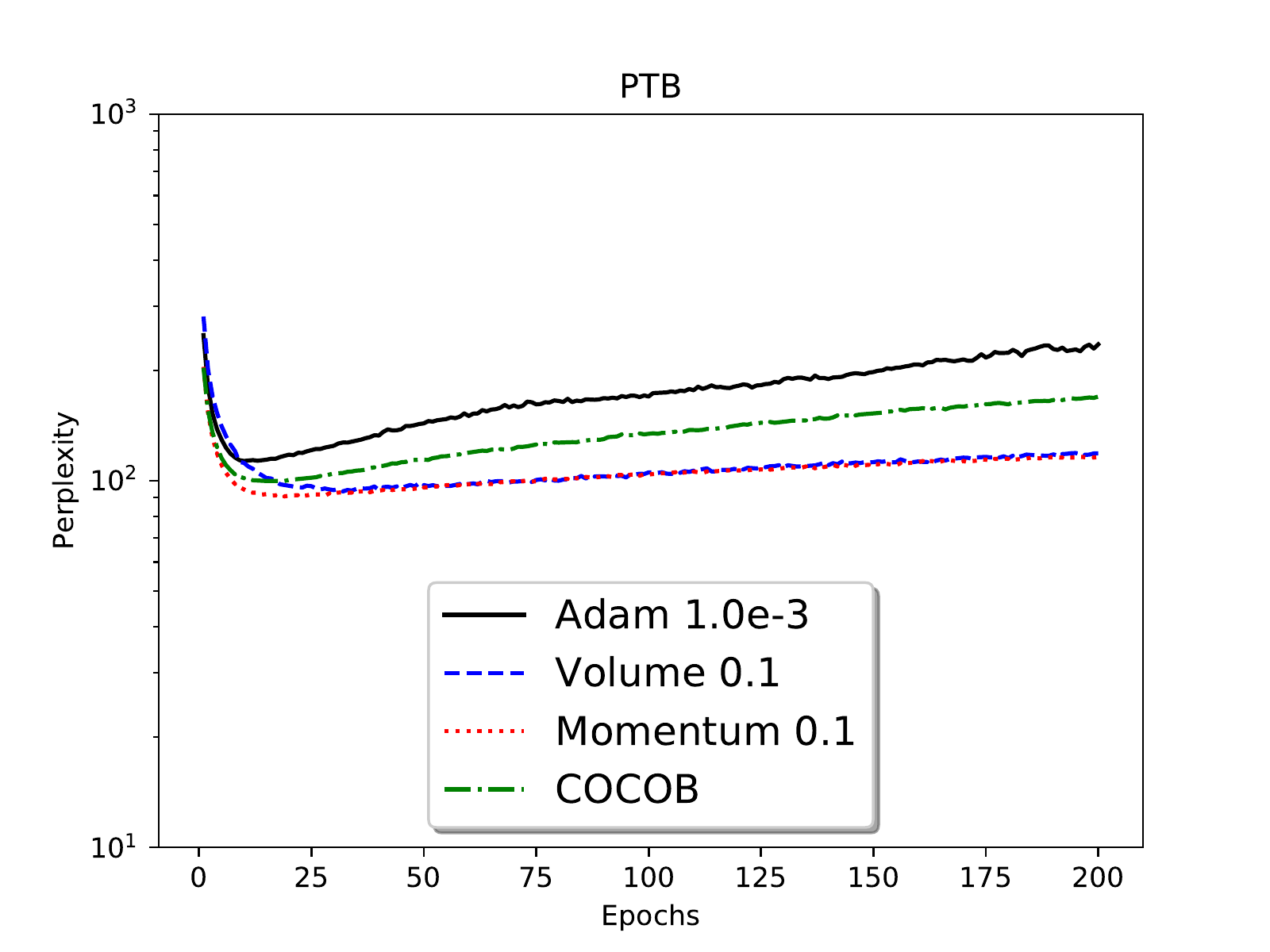}
\end{center}
\caption{PTB.}
\label{f6}
\end{figure}

Here we use the Penn Tree Bank dataset \cite{ptb} to train a Recurrent
Neural Network. The goal is to predict next words in a text given a history of previous words.
We use the code in the TensorFlow tutorial \cite{TFptb}. In Figure~\ref{f6} we plot
perplexity as a function of epochs for the train set and for the test set.

\subsection{Retraining Inception V3}
\begin{figure}[h!]
\begin{center}
\includegraphics[width=0.5\linewidth, height=0.53\linewidth]{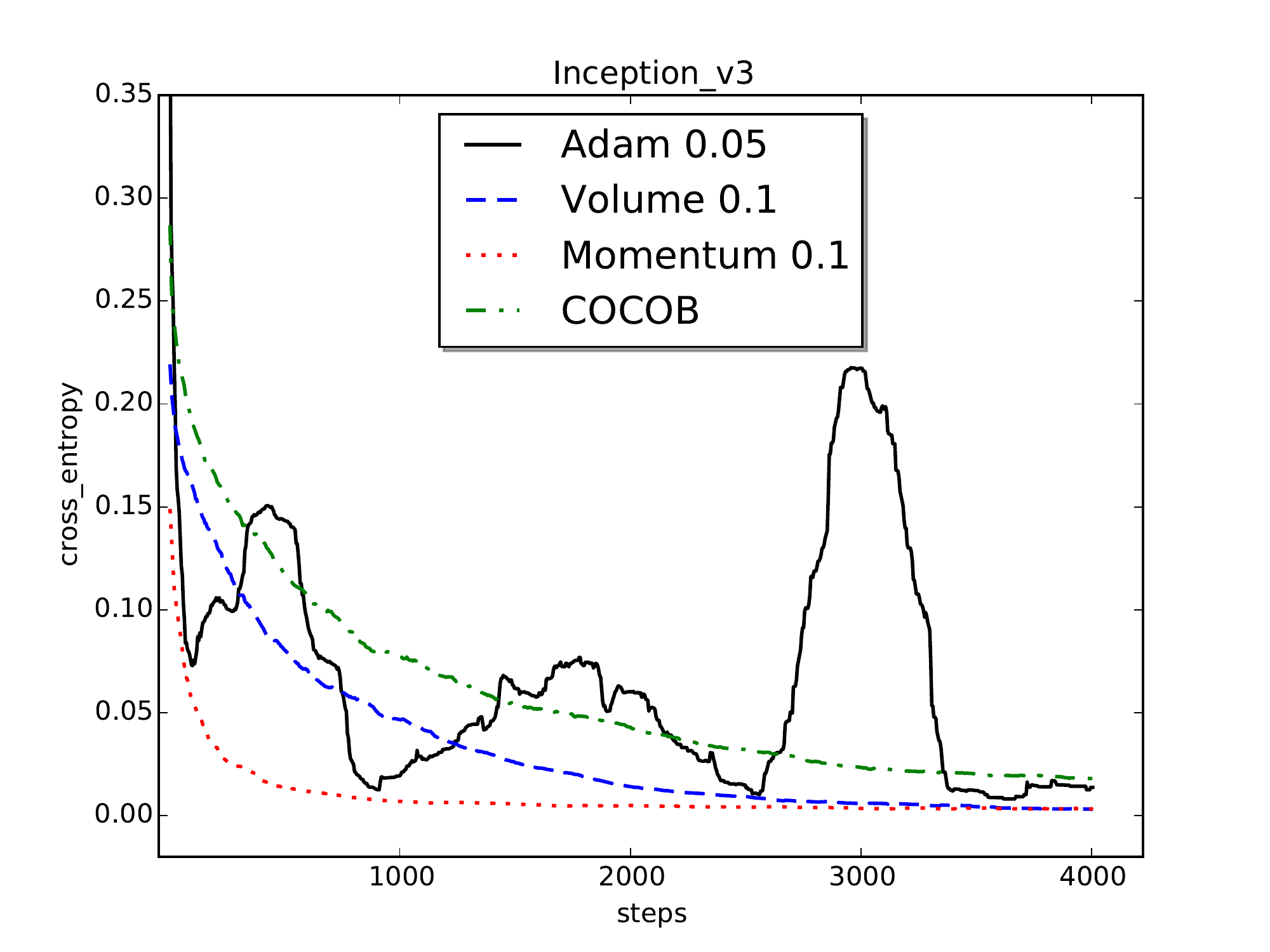}\includegraphics[width=0.5\linewidth, height=0.53\linewidth]{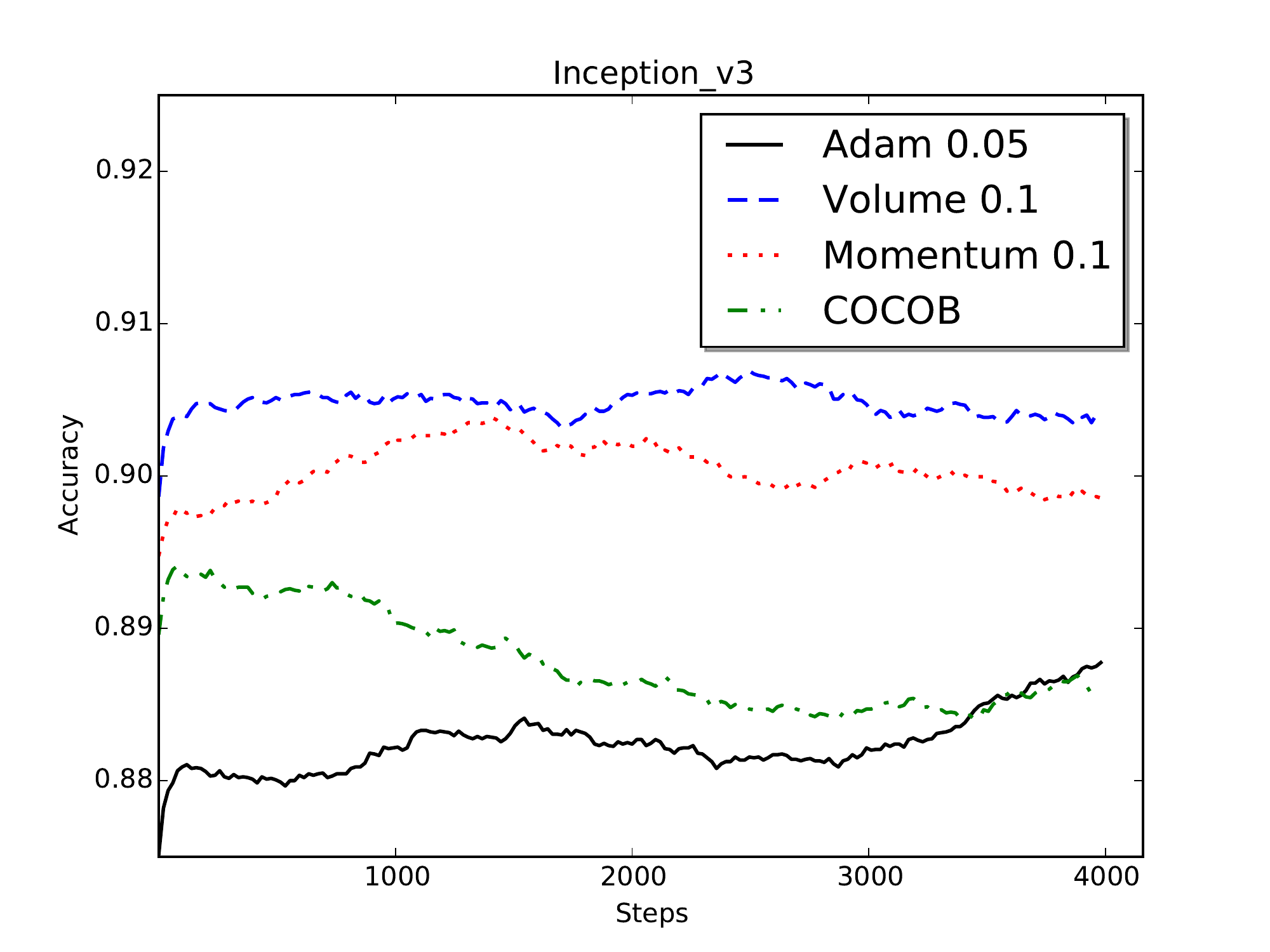}
\end{center}
\caption{Inception V3.}
\label{f7}
\end{figure}

Here we use a model trained on the Imagenet dataset \cite{imagenet} to classify images of flowers. We use the code from TensorFlow for Poets \cite{TFpoets} to download a pre-trained model, add a new final layer, and train that layer on a set of flower photos. Here we use an Inception V3 model, see \cite{V3}.
In Figure~\ref{f7} we plot cross entropy for the training set, and accuracy for
the evaluation set as a function of the number of steps.

\subsection{Retraining Mobilenet\_0.25\_160}
\begin{figure}[h!]
\begin{center}
\includegraphics[width=0.5\linewidth, height=0.5\linewidth]{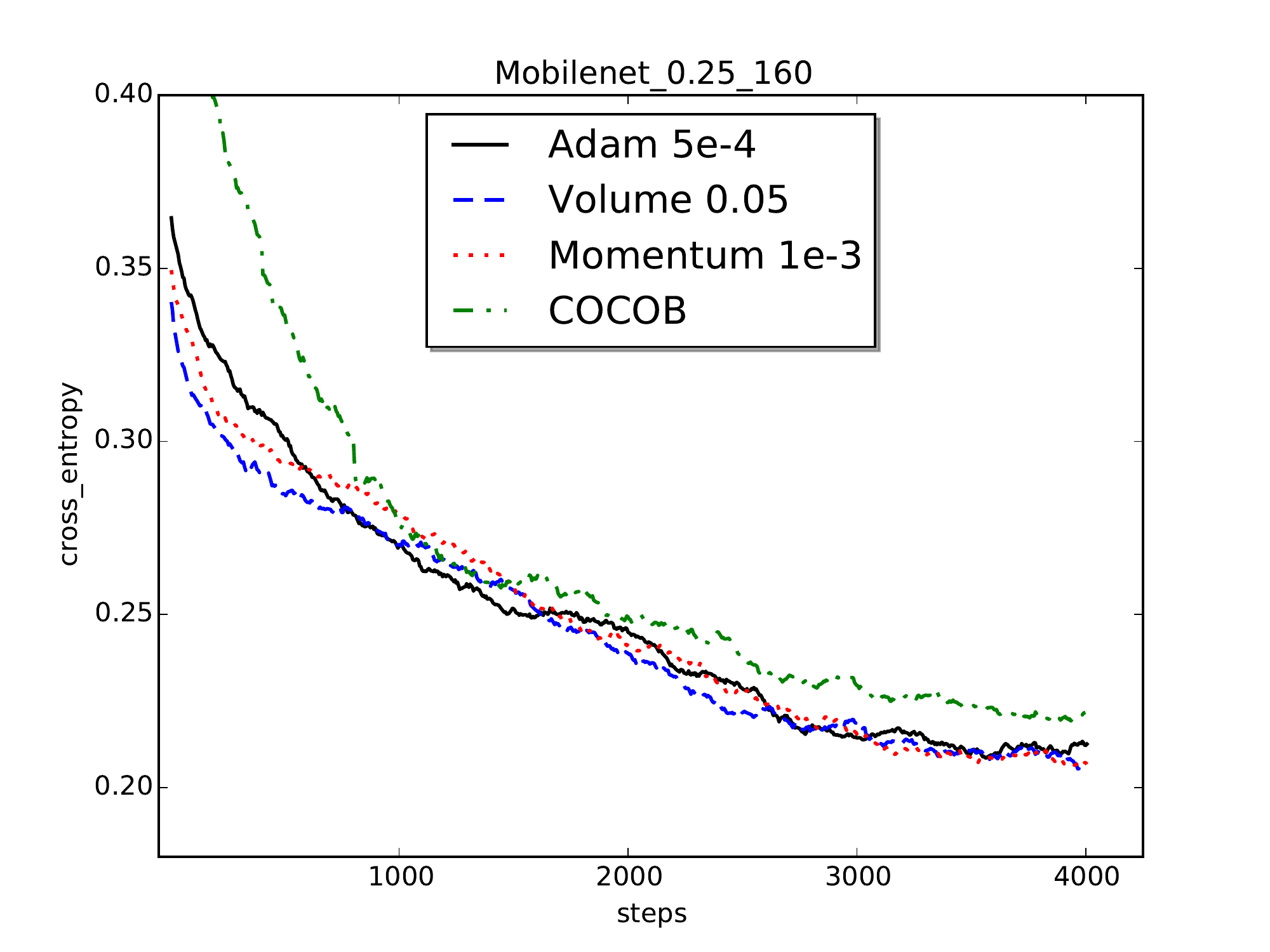}\includegraphics[width=0.5\linewidth, height=0.5\linewidth]{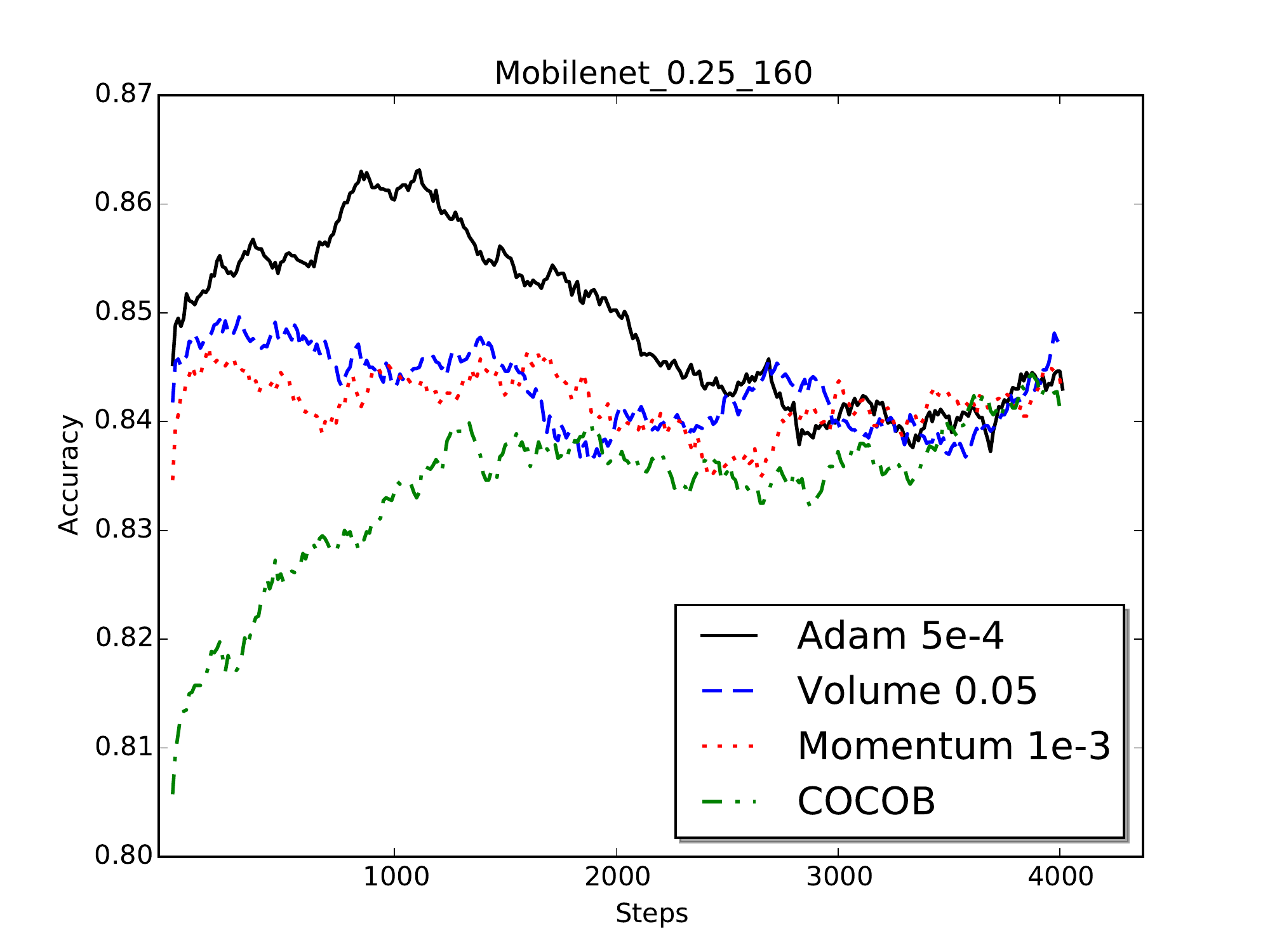}
\end{center}
\caption{Mobilenet\_0.25\_160.}
\label{f8}
\end{figure}

This time we retrain using a Mobilenet model \cite{mobilenet}. The relative size
of the model is
0.25 and the image size is 160. In Figure~\ref{f8} we plot cross entropy for the training set, and accuracy for
the evaluation set as a function of the number of steps.

\subsection{Retraining Mobilenet\_0.25\_128}
\begin{figure}[h!]
\begin{center}
\includegraphics[width=0.5\linewidth, height=0.5\linewidth]{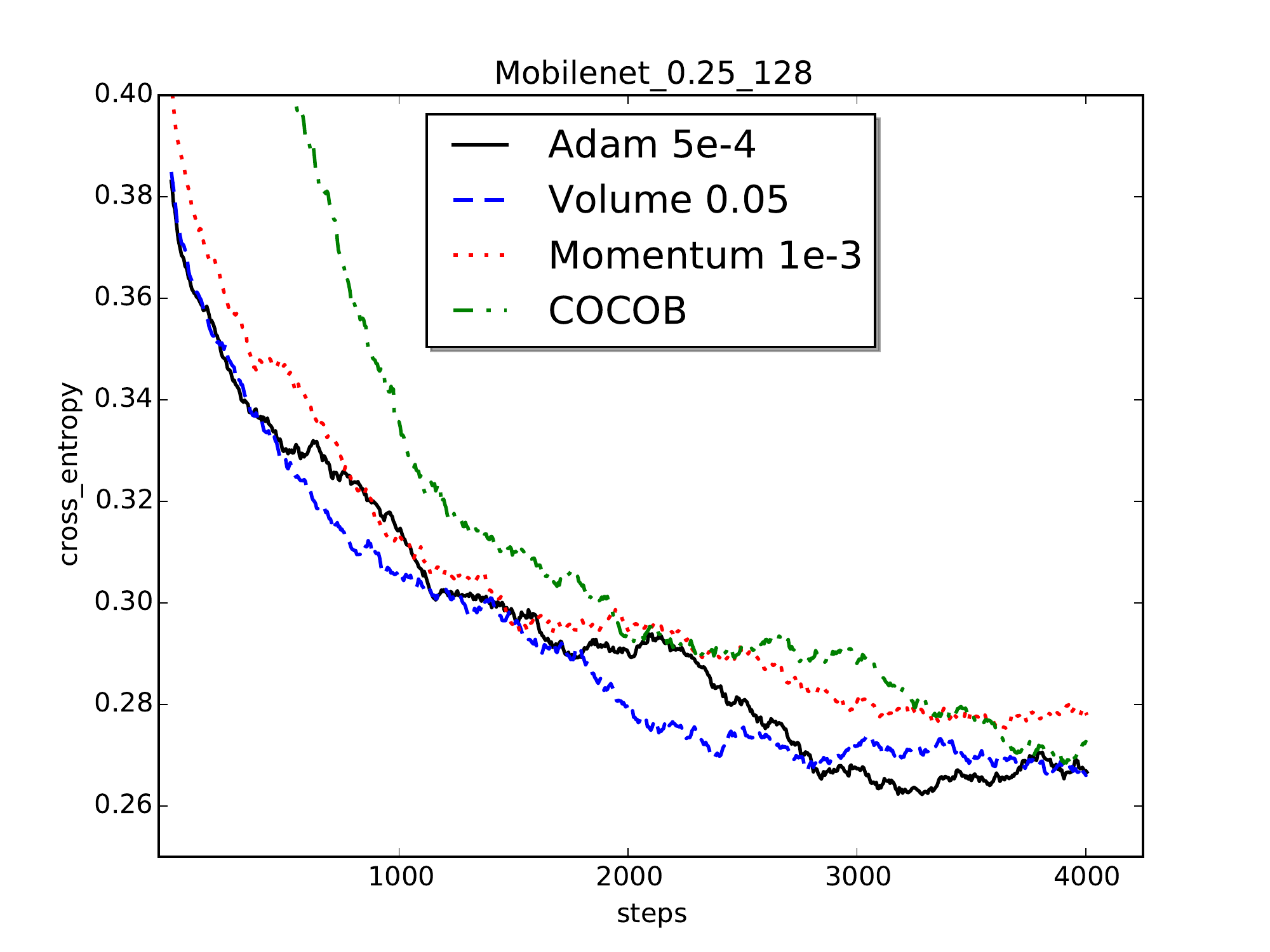}\includegraphics[width=0.5\linewidth, height=0.5\linewidth]{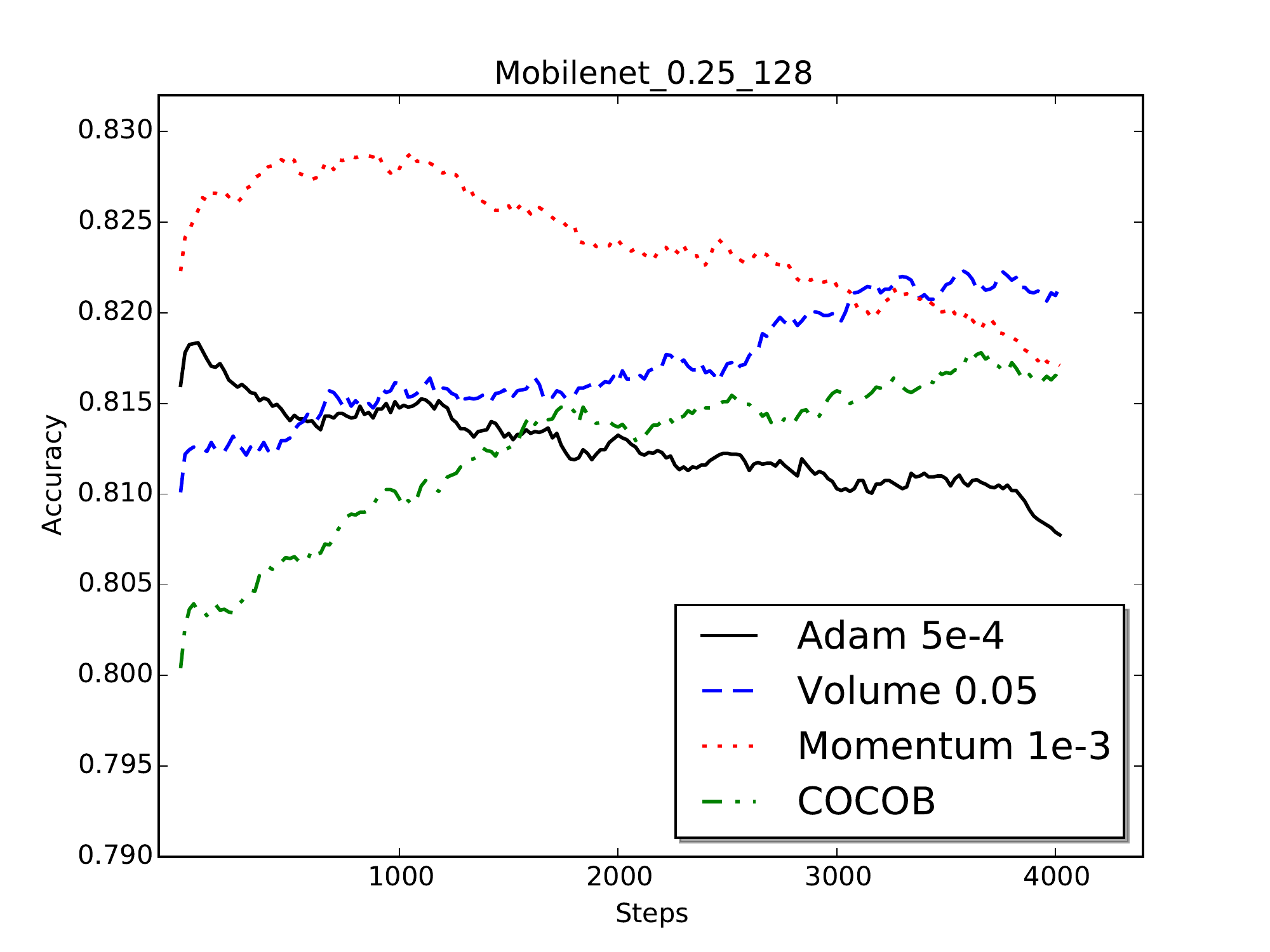}
\end{center}
\caption{Mobilenet\_0.25\_128.}
\label{f9}
\end{figure}

Here we retrained using a Mobilenet model \cite{mobilenet}. The relative size
of the model is
0.25 and the image size is 128. In Figure~\ref{f9} we plot cross entropy for the training set, and accuracy for
the evaluation set as a function of the number of steps.

\subsection{Simple Audio Recognition}
\begin{figure}[ht]
\begin{center}
\includegraphics[width=0.5\linewidth, height=0.48\linewidth]{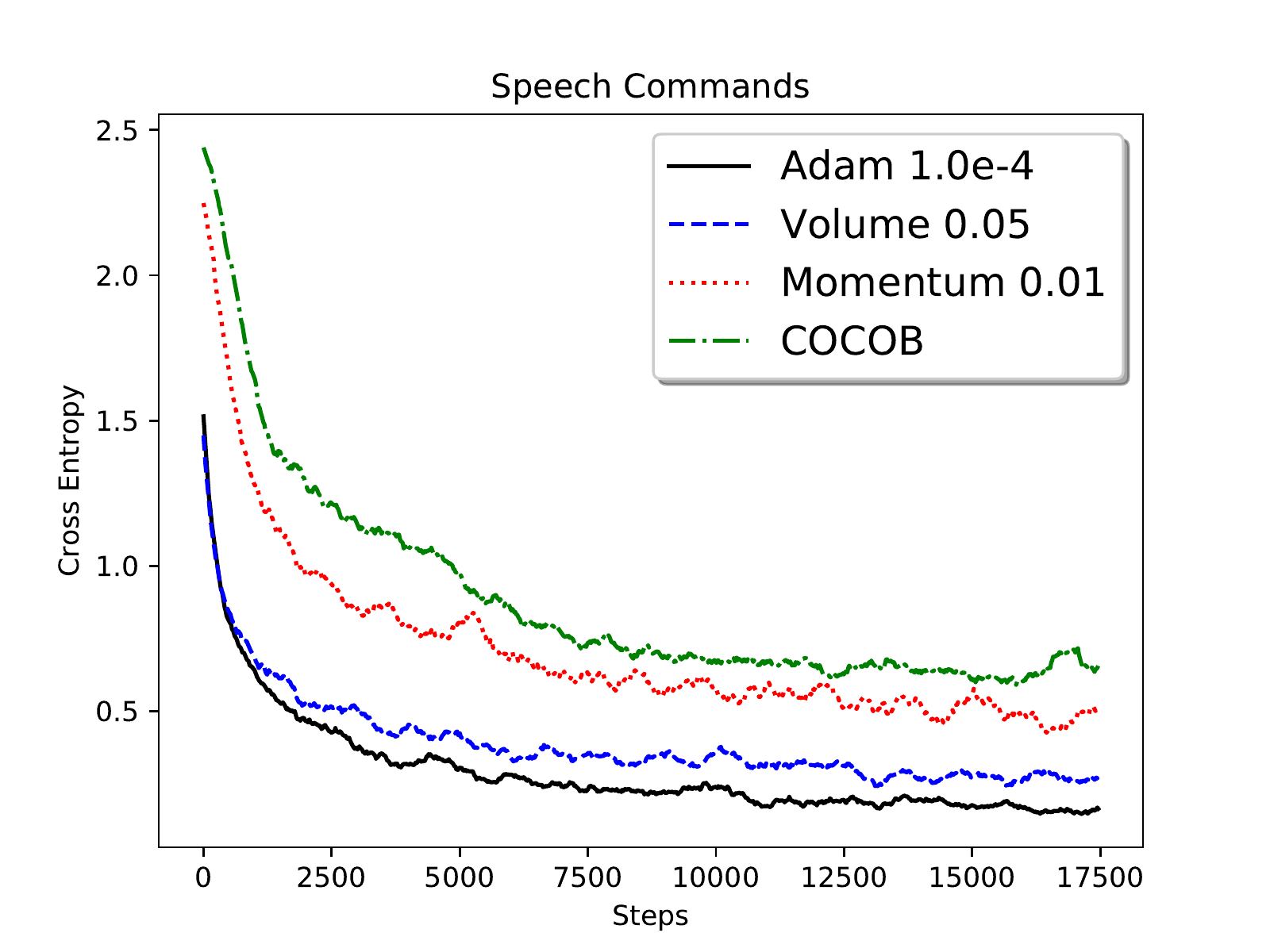}\includegraphics[width=0.5\linewidth, height=0.48\linewidth]{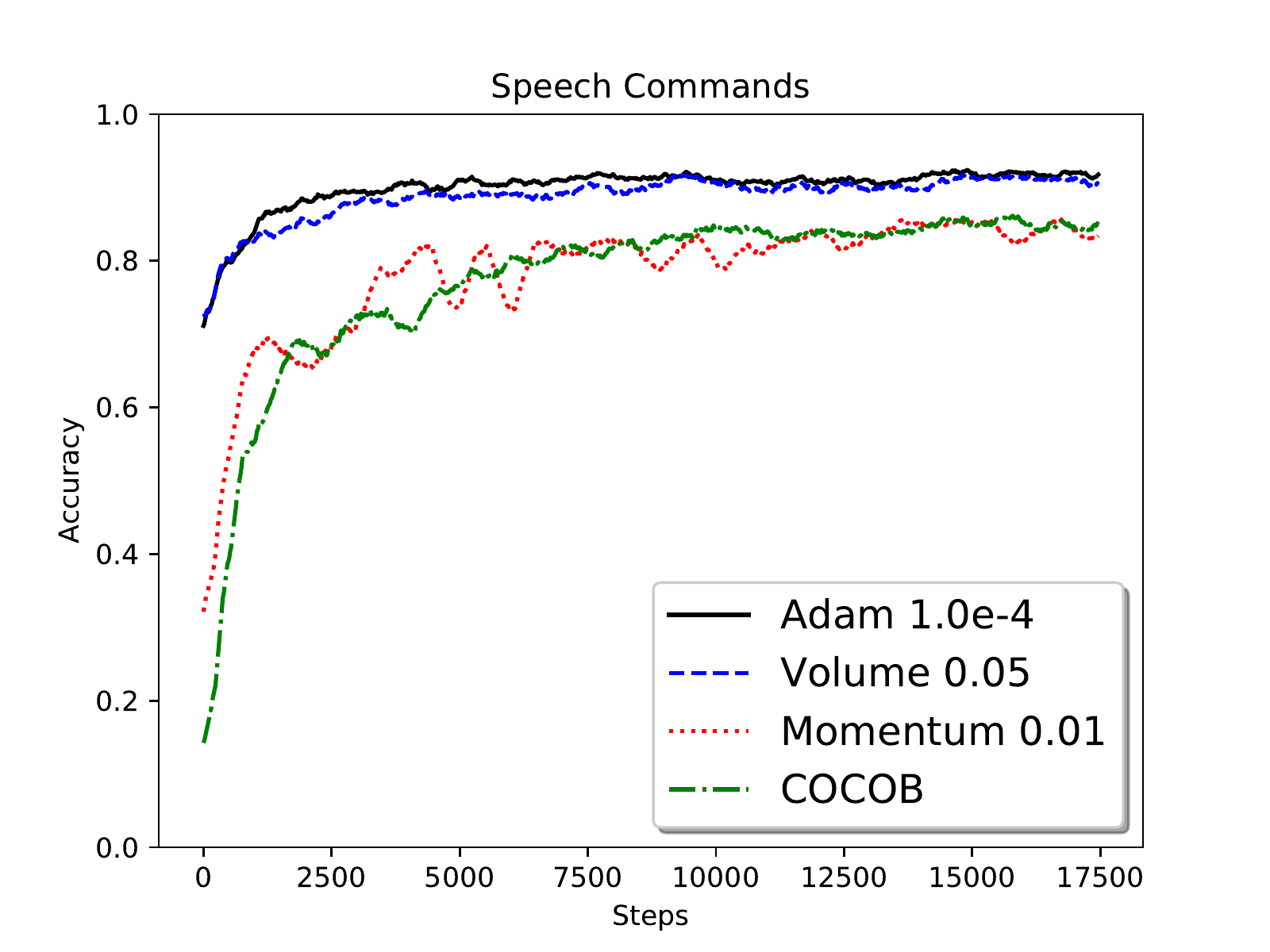}
\end{center}
\caption{Audio Recognition}
\label{f10}
\end{figure}

Here we use a basic speech recognition network that recognizes ten different words
from the TensorFlow tutorial \cite{speech}. In Figure~\ref{f10} we plot cross entropy for the training set, and accuracy for
the evaluation set as a function of the number of steps.

\section{Final Remarks}
We have presented a comparison of several algorithms for training
neural networks. In Table~\ref{sample-table} we give labels from
1 to 4 to
display the order
from best to worst of the different algorithms on each dataset.
If for one dataset several algorithms have been best, we give the
label 1 to each of them. The Volume Algorithm received the label 1
five times and the label 2 two times.

\begin{table}[h]
  \caption{Comparison of training algorithms}
  \label{sample-table}
  \centering
  \scalebox{0.92} {
  \begin{tabular}{lccccccc}
   \toprule
       & MNIST     & CIFAR10 & PTB & Inc V3 &  Mb\_.25\_160 & Mb\_.25\_128
       &Audio\\
    \midrule
    Adam& 2  & 3  & 2 & 3 & 1 & 1 &1\\
    Momentum   & 3  & 2 &  4 & 1 & 1 & 4 & 3\\
    COCOB   & 4   & 4   & 1 & 4 & 4 & 1 & 4\\
    Volume & 1 &1 & 2 & 1 & 1 & 1 & 2\\
    \bottomrule
  \end{tabular}
  }
\end{table}

\begin{figure}[ht]
\begin{center}
\includegraphics[width=7cm, height=4.cm]{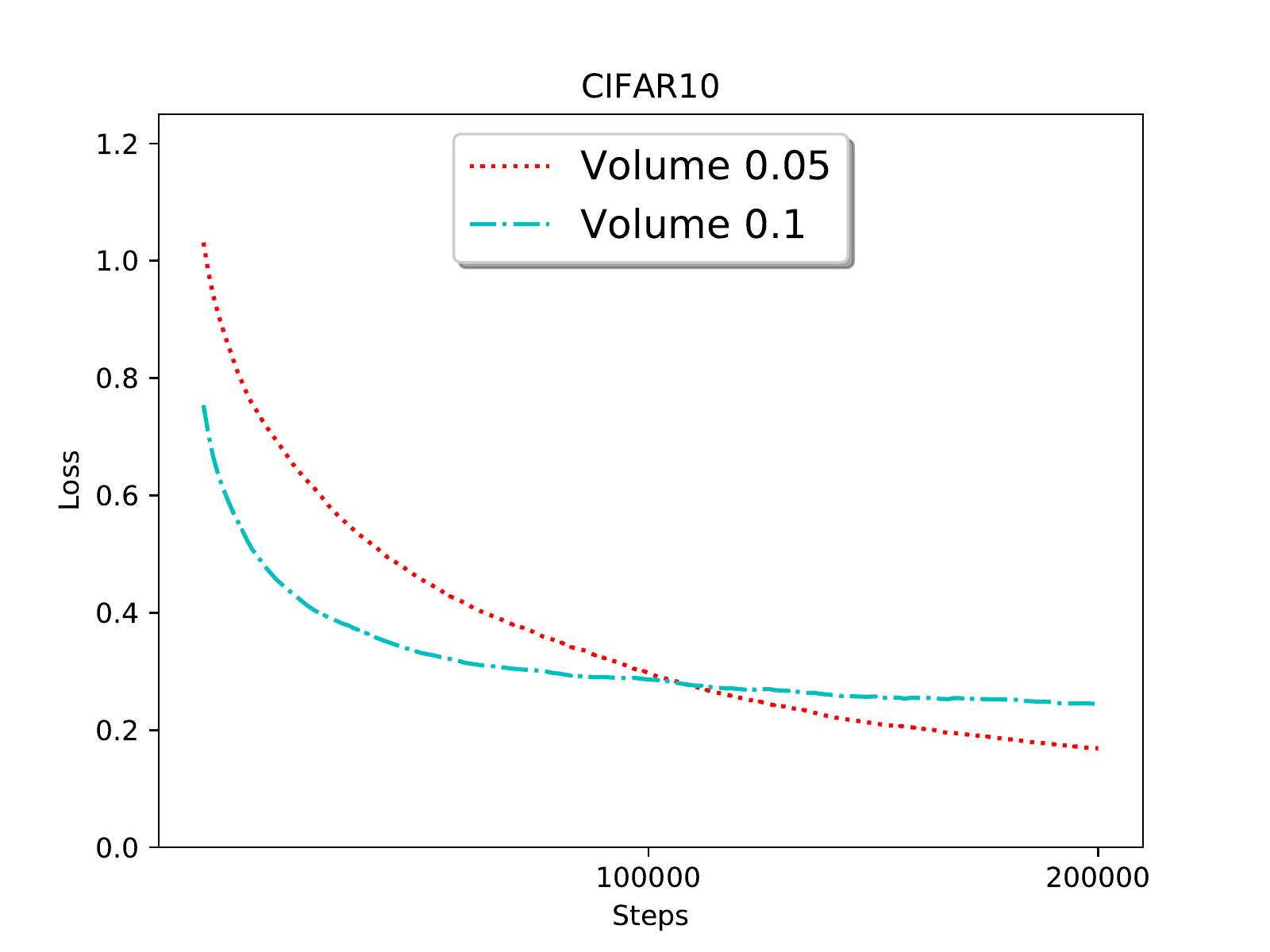}
\end{center}
\caption{Different step-sizes for the Volume Algorithm.}
\label{f11}
\end{figure}

The choice of the direction in our algorithm is similar to the one
for Adam, however the choice of the step-size makes a significant
difference, see Figure~\ref{f7} for instance. 
In our case we input an initial value $s$ and allow the step-size
to vary in the interval $[ 0.2s, 2s ]$. In Figure~\ref{f11} we show
the behavior of the Volume Algorithm for two different values of $s$,
namely the values $0.05$ and $0.1$. For the larger
value the function decreases faster in the earlier steps, and tends
to stabilize later. For the smaller value the decrease is slower
in the early steps, but it does not stabilize until later. This suggests
a strategy of using one interval for the early steps, and later switch
to an interval containing smaller values for the later steps.
The design of such strategy will be the subject of future research.

\section*{acknowledgments}
We are grateful Mark Wegman for his helpful comments on the preparation
of this paper.

\bibliographystyle{siam}

\bibliography{vol}

\end{document}